\documentclass[11pt]{article}
\pdfoutput=1
\usepackage{amsfonts, amssymb, amsmath, graphicx, epsfig, lscape, subfigure}
\usepackage{ifpdf}

\def\inst#1{$^{#1}$}

\newtheorem{theorem}{Theorem}[section]

\newtheorem{proposition}[theorem]{Proposition}

\newenvironment{proof}{\noindent\textbf{Proof:}}{\mbox{}\nobreak\hfill\qed}
  {\hspace*{\fill}$\Box$\par\vspace{4mm}}
  {\hspace*{\fill}$\Box$\par\vspace{4mm}}
  {\hspace*{\fill}\par\vspace{4mm}}

\newcommand{\qed}{$\Box$}
\newcommand{\UBu}{UBm_{11}}
\newcommand{\UBz}{UBm_{10}}

\newcommand{\DHn}{D^H_G(n_1)}

\newcommand{\NP}{{\it NP}}

\begin {document}

\title{Structural Bounds on the Dyadic Effect}

\author{%
Matteo Cinelli\inst{1} \and Giovanna Ferraro\inst{1} \and Antonio
Iovanella\inst{1}%
}

\date{}

\maketitle

\begin{center}
{\footnotesize 
\inst{1} Department of Enterprise Engineering\ \

University of Rome ``Tor Vergata''\\
Via del Politecnico, 1 - 00133 Rome, Italy.\\
\texttt{matteo.cinelli@uniroma2.it\\
giovanna.ferraro@uniroma2.it\\
antonio.iovanella@uniroma2.it\\}}
\end{center}

\begin{abstract}

In this paper we consider the dyadic effect introduced in complex networks when nodes are distinguished by a binary characteristic. Under these circumstances two independent parameters, namely dyadicity and heterophilicity, are able to measure how much the assigned characteristic affects the network topology. All possible configurations can be represented in a phase diagram lying in a two-dimensional space that represents the feasible region of the dyadic effect, which is bound by two upper bounds on dyadicity and heterophilicity. Using some network's structural arguments, we are able to improve such upper bounds and introduce two new lower bounds, providing a reduction of the feasible region of the dyadic effect as well as constraining dyadicity and heterophilicity within a specific range.
Some computational experiences show the bounds' effectiveness and their usefulness with regards to different classes of networks.
\vspace{0.5cm}

\end{abstract}

\section{Introduction}
\label{Intro}

Complex systems modelled as networks exhibit global structures that are commonly affected by the characteristics of their founding elements. Indeed, the properties assigned to these elements often correlate with the architecture of the observed systems (\cite{BBV}, \cite{NWB}). 
This paper is devoted to the situations where nodes themselves have peculiar properties that carry significant information regarding their role within the network topology. In literature (\cite{dealme}, \cite{MSC}), the tendency of nodes to link with others that are similar to themselves is a phenomenon called {\it homophily}, which affects the dyadic similarities between nodes and creates correlated patterns among neighbours. The nodes' tendency to connect with each other also relates to the concept of assortative mixing~\cite{New2} which describes correlations because of some nodes' properties. 

Park and Barab\'asi~\cite{PB} noted that, when nodes in a network fit within two distinct groups according to their characteristics, two different parameters, namely dyadicity and heterophilicity, are required to identify the relations between the network topology and nodes' features.

The dyadic effect has been considered in order to assess the functional role of nodes within biological networks such as, for instance, in gene-gene interaction in statistical epistasis networks~\cite{HAKM}, in phenome-genome networks~\cite{JDC} in disease-phenotype network~\cite{JLJZ} and in protein-protein interaction networks~\cite{ZZJ} where numerous characteristics are studied to evaluate genetic interactions. Nodes' characteristics are investigated also in inter-organizational innovation networks (\cite{FI2015}, \cite{FIP}) where partnerships agreements of technological transfer among countries are related to innovation indices.

The methodology presented in~\cite{PB} is able to study exhaustively all the configurations of a binary characteristic on a network using a phase diagram, which lies into a two dimensional space constrained by certain network related bounds. However, current bounds are computed a-priori considering either particular networks arguments or the number of featured nodes, which results in a space that is often much larger than necessary. Another important issue is related to the phase diagram computational complexity which grows exponentially in the number of nodes, implying some limitations in real applications. The literature has attempted to overcome these difficulties by using heuristics or statistical methods. For instance, in~\cite{PB} it is reported a heuristic method able to identify extremal configurations; in~\cite{HAKM}, ~\cite{JDC}, ~\cite{JLJZ}, and~\cite{ZZJ} statistical methods are used to infer the existence of some configurations in computational biology, while in~\cite{BPM} entropy-based measures are used to globally assess the relevance of nodes' characteristics. 

The contribution of this paper lies in the improvement of current upper bounds through considerations related to structural arguments surrounding a given network. The reduction of the two dimensional space is performed by not only considering the two upper bounds but by also introducing two lower bounds. We also present the analytical reasonings and we test their behaviour on different classes of networks. 

New bounds' foundations are rooted in the degree sequence which can be easily extracted from any network. Although the results we obtain in the space reduction can be valuable, we make use of the straightforward relationships that allow us to compute such bounds independently of the network size. 

We provide a reduction of the dyadic effect's feasible region, whic can be used in all the applications where the degree of correlation of the nodes' characteristics with the network topology is evaluated considering empirical arguments instead of computing the phase diagram such as, for instance, in~\cite{CFI} and~\cite{HAKM}. Moreover, by introducing the four bounds we are limiting the values of dyadicity and heterophilicity to lie in a range that can be easily computed.

The paper is organized as follows: Section~\ref{probset} gives the problem setting; Section~\ref{ULB} shows the upper and lower bounds; Section~\ref{EE} gives the computational analysis; Section~\ref{conclusion} presents the conclusions.

\section{Problem settings}\label{probset}
\subsection{Theoretical Background}
\label{ThBg}

The classical mathematical abstraction of a network is a graph $G$. Let $G$ = $(V,E)$ be a graph composed of a set $V$ of $N$ nodes and a set $E$ of $M$ edges that defines the relationship among these nodes. Herein $G$ is considered undirected, unweighted, connected and simple, i.e. loops and multiple edges are not allowed. 

The degree $d_i$ of a node $i \in V$ is defined as the number of edges in $E$ incident to $i$. The nodes' degrees listed in an non-increasing order are referred to as \textit{degree sequence} $D_G$ and, as we recall, for every connected graph holds the Degree-Sum Formula or \textit{Handshaking Lemma},  $\sum_{i=1}^{N}{d_i} = 2M$. A \textit{graphic sequence} is defined as a list of nonnegative numbers which is the degree sequence of certain simple graphs. A graph $G$ with degree sequence $D_G$ is called a \textit{realization} of $D_G$.

A generic list $L$ of nonnegative numbers is not necessarily a graphic sequence. Indeed, a necessary and sufficient condition is that $\sum_{i=1}^{N}{d_i}$ is even and  $\sum_{i=1}^{N}{d_i} \leq k(k-1)+ \sum_{i=k+1}^n \min\{k, di\}$, $1 \leq k \leq N$ ~\cite{EG}. The problem of discovering if $L$  is a graphic sequence is called \textit{Graph Realization Problem}~\cite{H1962, H1955}.

Given an integer $n \leq N$, we consider the subsequence of the first $n$ elements of $D_G$ calling it $D_G^H(n) \subseteq D_G$ as the head of $D_G$, and the subsequence of the last $n$ elements of $D_G$ calling it $D_G^T(n) \subseteq D_G$ as the tail of $D_G$.

A \textit{clique} $K_n$ is a complete subgraph of $G$ of dimension $n$, i.e. a subgraph of $n$ mutual interconnected nodes. The problem of finding the clique of highest cardinality in $G$ is a well known \NP-problem~\cite{GJ1976}. Since the degree of nodes in a clique of cardinality $n$ is at least $n-1$, a necessary condition for the existence of such clique is that $D_G$ contains at least $n$ nodes of degree $d_i \geq n -1$, otherwise the graph $G$ can not contain such clique.

A \textit{star} $S_n$ is a subgraph of $G$ of dimension $n$ showing one node with degree $n-1$ and the others $n -1$ having degree $1$.

Traditionally, in graph theory, edges have two endpoints since, by definition, they represent a reciprocal relationship between two nodes. In the literature of graph realisation problem, as well as that of dynamic graphs, it is possible to admit half-edges anchored at one node of a degree sequence while the other endpoint is free. This particular object is called \textit{stub}~\cite{NWS}; when two stubs of two distinct nodes connect, then a classical edge is realized~\cite{KTEMS}.

\subsection{Nodes' characteristics and dyadic effect}
\label{nodechar}

\begin{figure}[htbp]
\begin{center}
\includegraphics[trim=4cm 24.3cm 4cm 3.5cm, clip = true, totalheight=0.07\textheight]{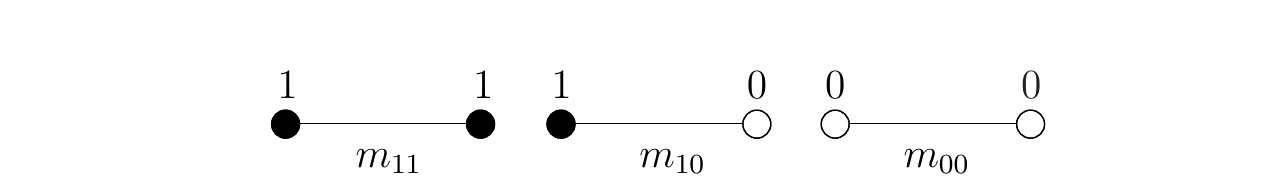}
\caption{Types of dyads.}
\label{bpd}
\end{center}
\end{figure}

Herein, we refer to a given characteristic $c_i$, which can assume the values $0$ or $1$, for each $i \in N$. Consequently, $N$ can be divided into two subsets: the set of $n_1$ nodes with characteristic $c_i = 1$, the set of $n_0$ nodes with characteristic $c_i = 0$; thus, $N = n_1 + n_0$. We distinguish three kinds of \textit{dyads}, i.e. edges and their two end nodes, in the network: $(1 -� 1)$, $(1 -� 0)$, and $(0 -� 0)$ as depicted in the Figure~\ref{bpd}.

We label the number of each dyad in the graph as $m_{11}$, $m_{10}$, $m_{00}$, respectively. Hence, $M = m_{11} + m_{10} + m_{00}$. We consider $m_{11}$ and $m_{10}$ as independent parameters that represent the dyads containing nodes with characteristic 1.

\begin{figure}[htbp]
\begin{center}
 \begin{minipage}[b]{6cm}
   \centering
   \includegraphics[scale=.8]{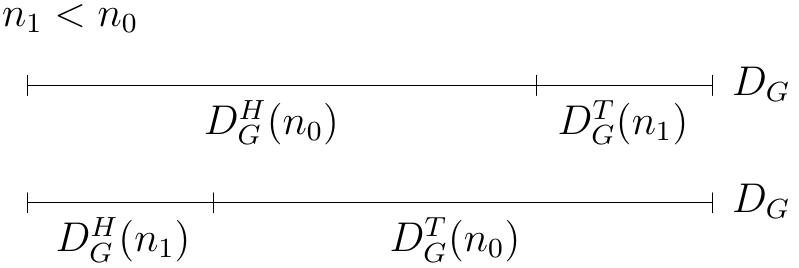}
   \caption{Two different partitions when $n_1 < n_0$.}\label{deg1}
 \end{minipage}
 \ \hspace{5mm} \hspace{5mm} 
 \begin{minipage}[b]{6cm}
  \centering
   \includegraphics[scale=.8]{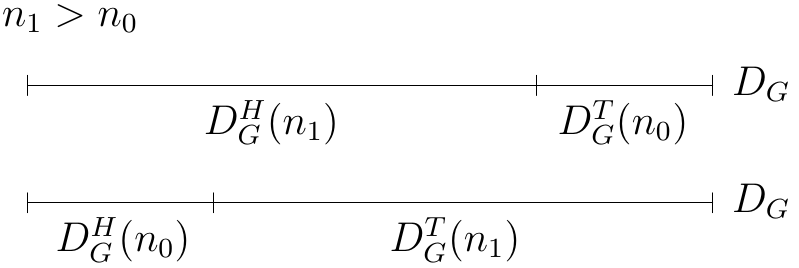}
   \caption{Two different partitions when $n_1 > n_0$.}\label{deg2}
 \end{minipage}
 \end{center}
\end{figure}

Let $D_G$ be the degree sequence of $G$, we can use $n_1$ and $n_0$ to define its heads $D_G^H(n_1)$ or $D_G^H(n_0)$ and the tails $D_G^T(n_1)$ or $D_G^T(n_0)$ such as $D_G = D_G^H(n_1) \cup D_G^T(n_0) = D_G^H(n_0) \cup D_G^T(n_1)$. These partitions of the degree sequence are given arbitrarily assigning the characteristic $c_i = 1$  to the $n_1$ nodes with the highest degree or to the $n_1$ nodes with lowest degree or viceversa. Such partitions are reported in Figures~\ref{deg1} and~\ref{deg2}, distinguishing the case in which $n_1 < n_0$ or $n_1 > n_0$. We make this construction in order to use it in Section~\ref{ULB}.

When nodes in a network fit within two distinct groups according to their characteristics, two different parameters are required to determine the existence of the relations between the network topology and the nodes' features~\cite{PB}. In many systems, the number of edges between nodes sharing a common characteristic is larger than expected if the characteristics are distributed randomly on the graph; this phenomenon is called the dyadic effect~\cite{WH}.
If a casual setting among the $N$ nodes is considered, where any node has an equal chance of having the characteristic 1, the values of $m_{11}$ and $m_{10}$ are~\cite{PB}:

\begin{equation}
\overline{m}_{11} = {n_1 \choose 2} \delta = \frac{n_1 (n_1 - 1)}{2} \delta
\label{m11}
\end{equation}

\begin{equation}
\overline{m}_{10} = {n_1 \choose 1} {n_0 \choose 1} \delta = n_1 (N - n_1) \delta
\label{m10}
\end{equation}

\noindent
where $\delta$ is the density and is equal to $\delta = 2M/N(N - 1)$. The relevant deviations of $m_{11}$ and $m_{10}$ from the expected values 
$\overline{m}_{11}$ and $\overline{m}_{10}$ denote that the characteristic $1$ is not randomly distributed~\cite{dealme, PB}. Such deviations can be calculated through the ratios of dyadicity $D$ and heterophilicity \textit{H} defined as: 
\begin{equation}
D = \frac{m_{11}}{\overline{m}_{11}}
\label{eq_D11}
\end{equation}

\begin{equation}
H = \frac{m_{10}}{\overline{m}_{10}}
\label{eq_H}
\end{equation}

If the characteristic is dyadic, $D > 1$, it means that nodes with the same characteristics tend to link more tightly among themselves than expected in a random configuration. Conversely when $D < 1$, the characteristic is anti-dyadic, indicating that similar nodes tend to connect less densely among themselves than expected in a random configuration. The characteristic is defined as heterophilic, with a value $H > 1$, highlighting that nodes with the same features have more connections to nodes with different characteristics than expected randomly. On the contrary, with a value $H < 1$, the characteristic is defined as heterophobic, meaning that nodes with certain characteristics have fewer links to nodes with different characteristics than expected randomly.

In~\cite{PB}, it is enstablished that $m_{11}$ and $m_{10}$ cannot assume arbitrary values, as there are indirect constraints due to the network structure. Indeed, $m_{11}$ cannot exceed 
\begin{equation}\label{ub11old}
\UBu = min( M, \binom { n_1 } {2} )
\end{equation}
\noindent
and $m_{10}$ cannot be larger than 
\begin{equation}\label{ub10old}
\UBz = min( M, n_1 n_0)
\end{equation}
\noindent
where \textit{UB} stands for upper bound.

\subsection{The phase diagram}
\label{phasdiag}

One instrument to investigate the correlation among the distribution of a given property $c$ and the underlying network structure is the phase diagram which, in general, describes the admissible configurations in the graph.

\begin{figure}
\begin{center}
 \begin{minipage}[b]{6cm}
   \centering
   \includegraphics[trim=7cm 19cm 3cm 4.3cm, clip=true, totalheight=0.23\textheight]{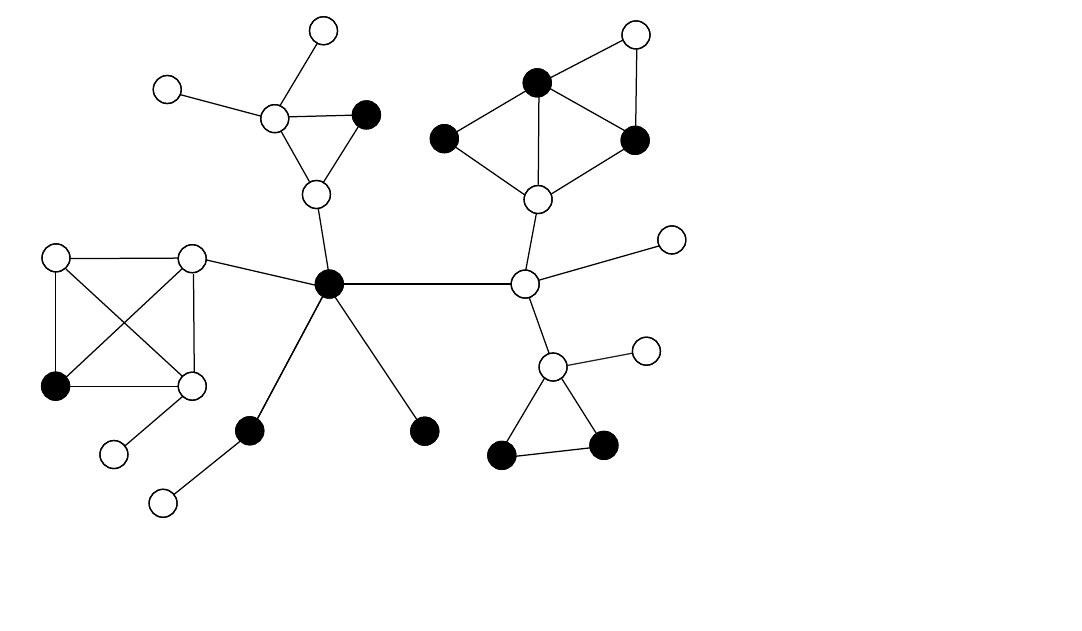}
   \caption{An example of graph with $N = 25$,  $M = 32$ and $n_1 = 10$.}\label{example}
 \end{minipage}
 \ \hspace{2mm} \hspace{2mm} 
 \begin{minipage}[b]{6cm}
  \centering
   \includegraphics[scale=.38]{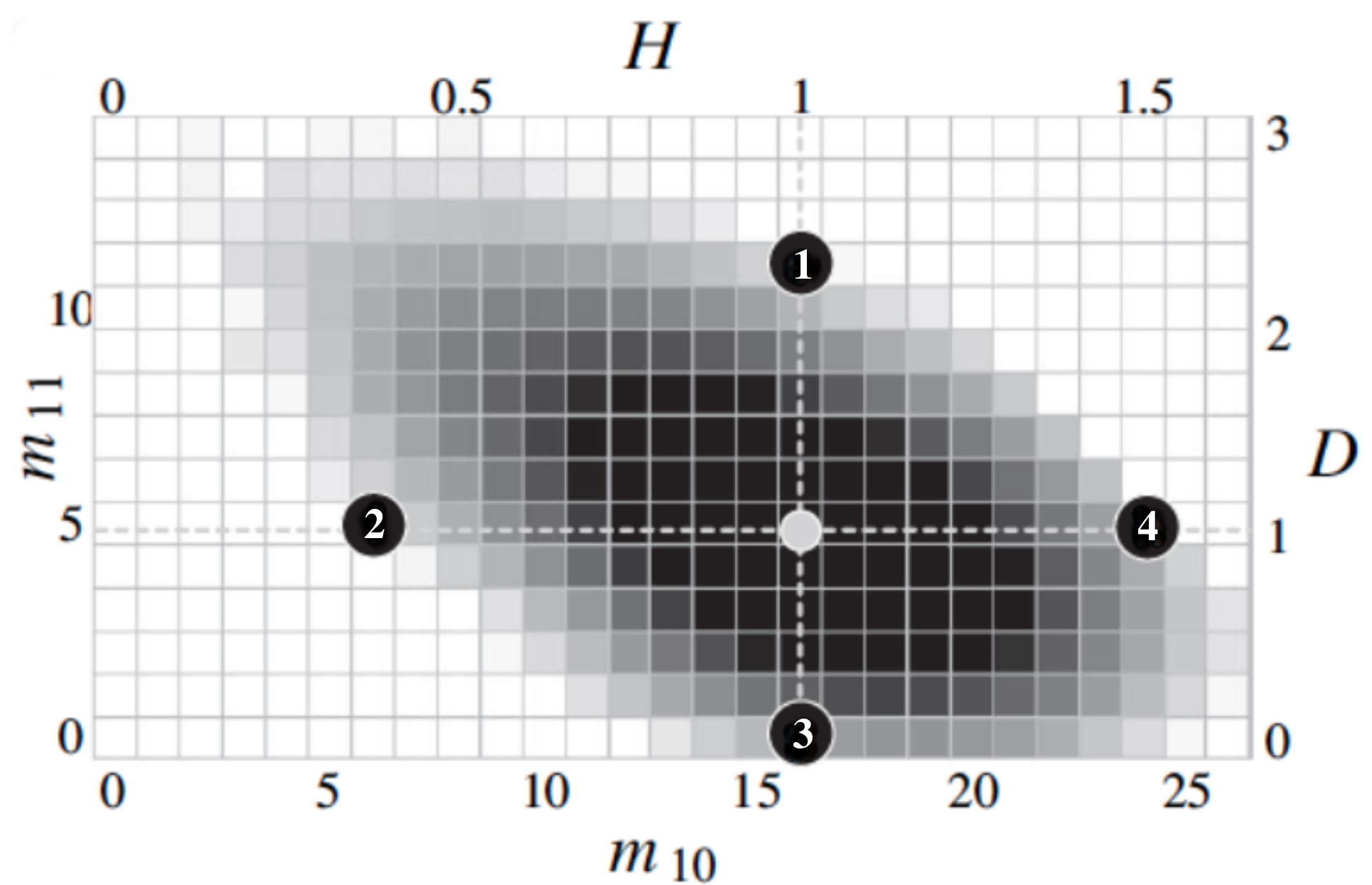}
   \caption{The phase diagram of possible values of ($m_{11}$, $m_{10}$).}\label{pdexample}
 \end{minipage}
 \end{center}
\end{figure}

\begin{figure}
\begin{center}
 \begin{minipage}[b]{3cm}
   \centering
   \includegraphics[trim=7cm 18.5cm 3cm 3cm, clip=true, totalheight=0.17\textheight]{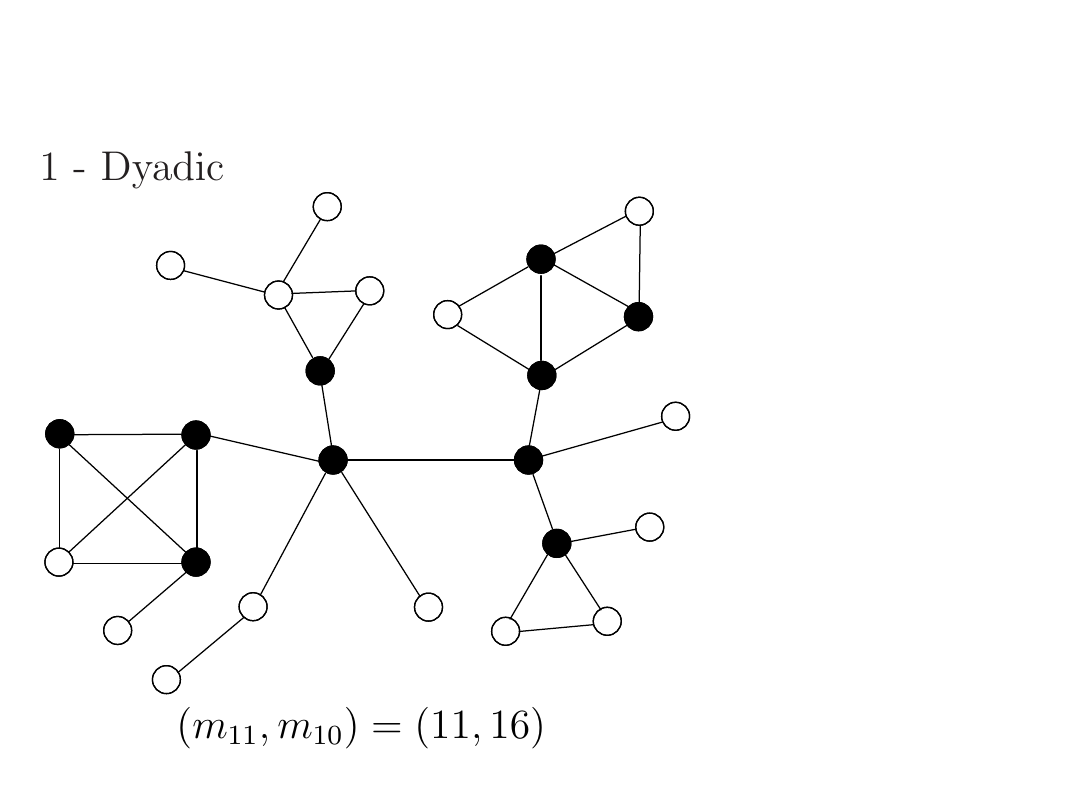}
 \end{minipage}
\ \hspace{1mm} \hspace{1mm} \
 \begin{minipage}[b]{3cm}
  \centering
   \includegraphics[trim=7cm 18.5cm 3cm 3cm, clip=true, totalheight=0.17\textheight]{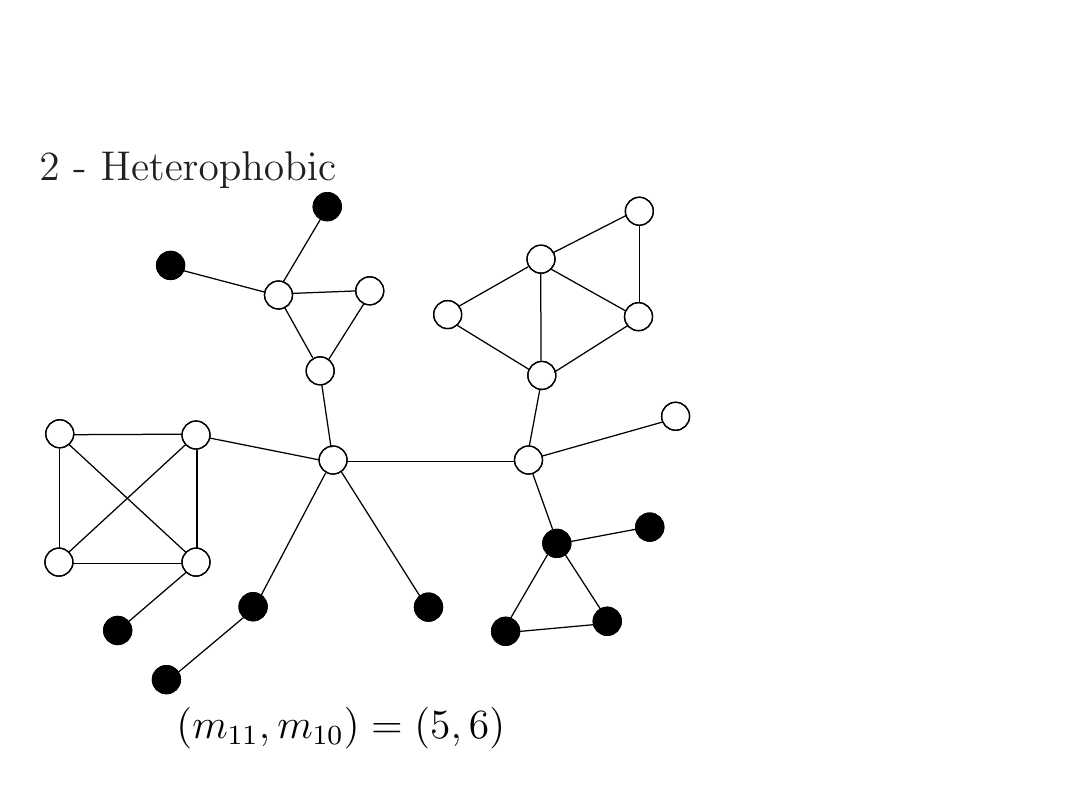}
 \end{minipage}
 \ \hspace{1mm} \hspace{1mm} \
 \begin{minipage}[b]{3cm}
  \centering
   \includegraphics[trim=7cm 18.5cm 3cm 3cm, clip=true, totalheight=0.17\textheight]{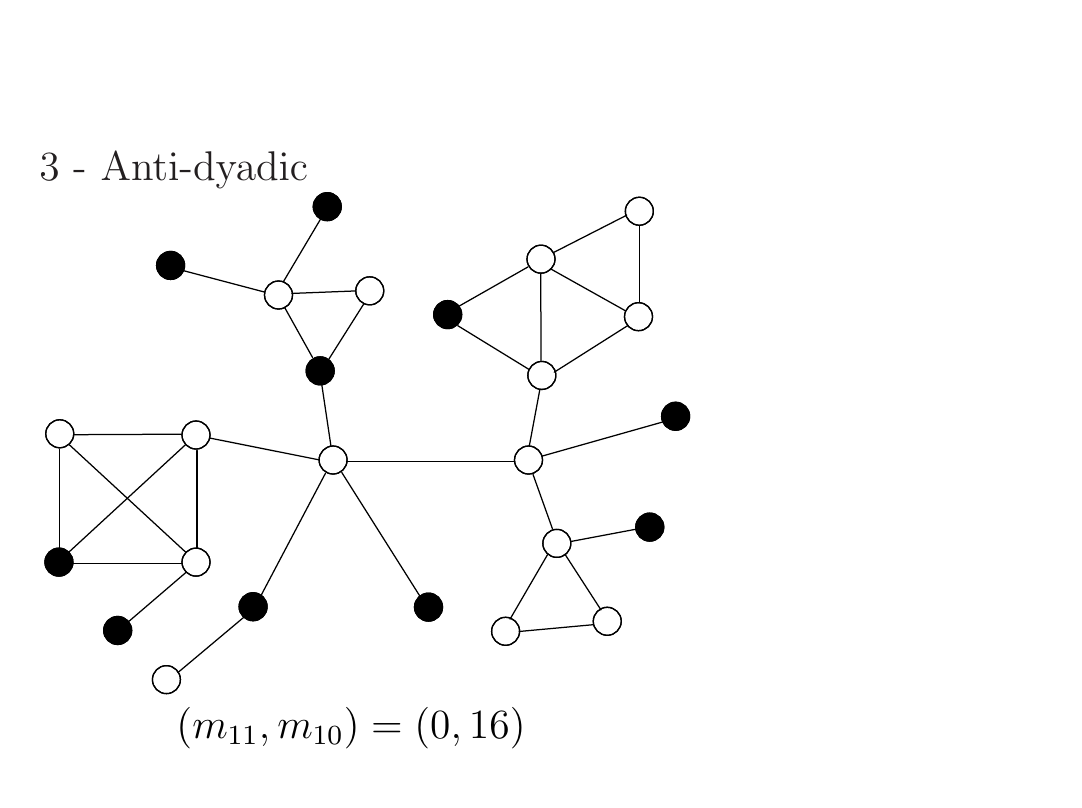}
 \end{minipage}
 \ \hspace{1mm} \hspace{1mm} \
  \begin{minipage}[b]{3cm}
  \centering
   \includegraphics[trim=7cm 18.5cm 3cm 3cm, clip=true, totalheight=0.17\textheight]{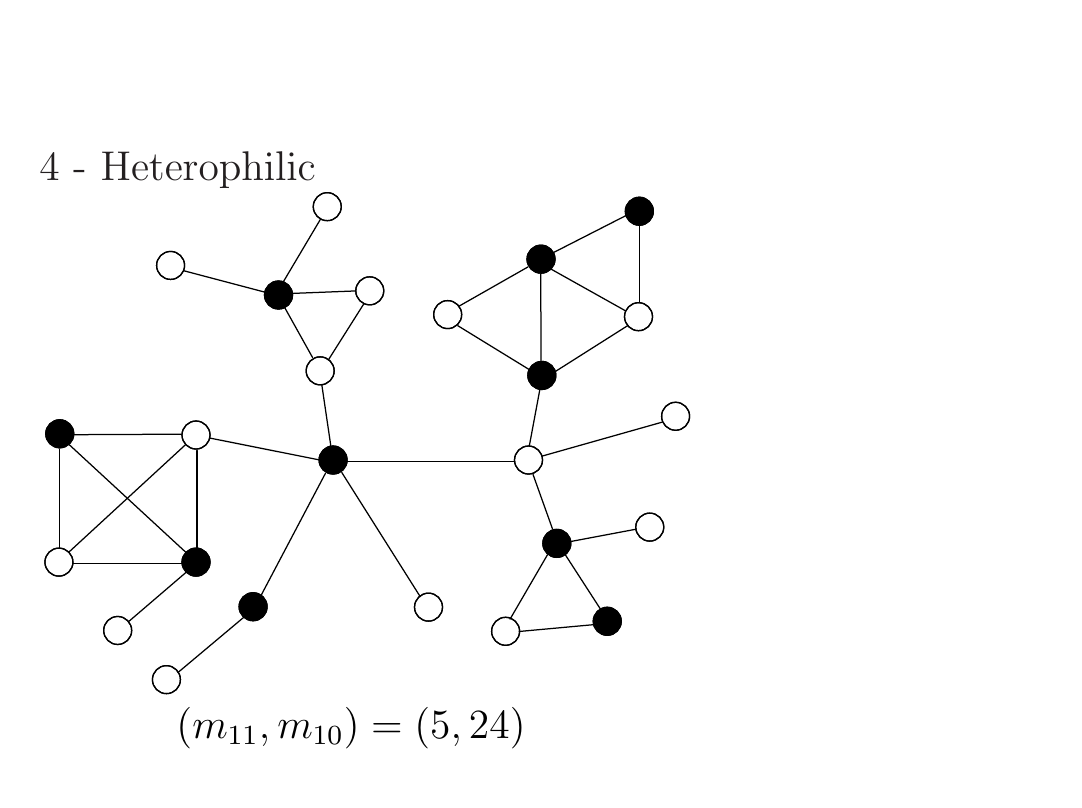}
 \end{minipage}
\caption{Configurations of four extreme points on the phase diagram.}
\label{fig_config}
\end{center}
\end{figure}

We consider, as an example, the graph shown in Figure~\ref{example} , which depicts a network with $25$ nodes and $32$ edges of which $n_1 = 10$ black nodes are randomly distributed (example and Figures~\ref{example}--~\ref{fig_config} are taken from~\cite{PB}).

The corresponding phase diagram in Figure~\ref{pdexample} describes the distribution of a random feature in the system. It should be noted that the reported phase diagram shows only a subarea.

The phase diagram presents all the admissible combinations of $m_{10}$ ($x$-coordinate) and $m_{11}$ ($y$-coordinate) and each corresponding square collects the number of assignments of $n_1$ nodes over the set $N$ for every fixed $m_{10}$ and $m_{11}$. There is a direct correspondence among the $m_{10}$ and $m_{11}$ axis and, respectively, $H$ and $D$, since the values are related through means of Eq.~\ref{eq_D11} and Eq.~\ref{eq_H}. Moreover, $m_{10}$ ranges from $0$ to $UBm_{10}$ and $m_{11}$ from $0$ to $UBm_{11}$. Correspondingly, $D$ ranges from $0$ to $D_{max} = UBm_{11}/ \overline{m}_{11}$ and $H$ ranges from $0$ to $H_{max} = UBm_{10}/\overline{m}_{10}$.

\begin{figure}[htbp]
\begin{center}
\includegraphics[scale=.6]{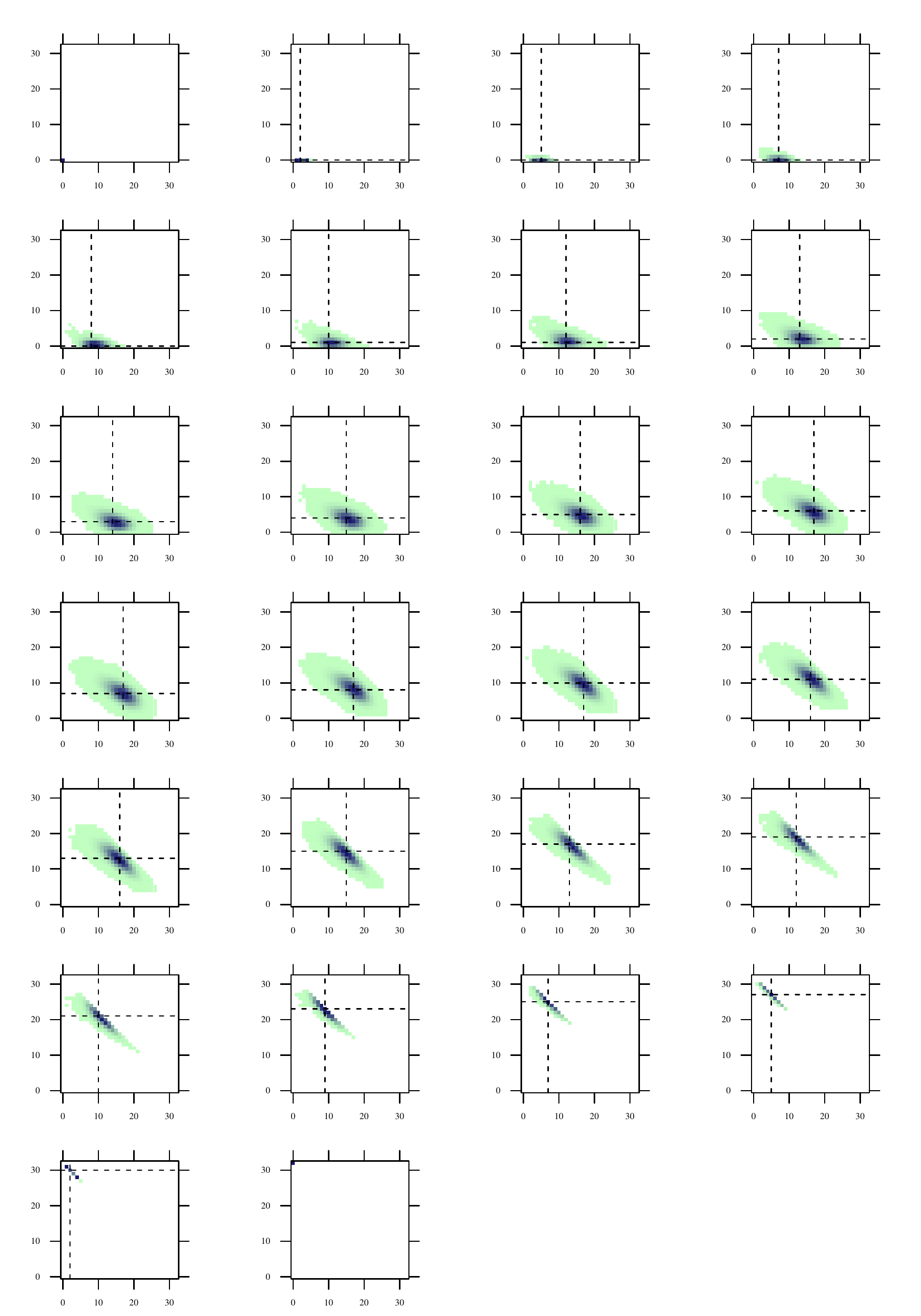}
\caption{The sequence of the phase diagrams of the network depicted in Figure~\ref{example}, starting from left to right and from the upper to the lower side of the figure. Dotted lines indicate the coordinates of $D = H = 1$.}
\label{pdsequence}
\end{center}
\end{figure}

For a given $m_{10}$ and $m_{11}$, each square has a darkness proportional to the degeneracy of the configuration and an open square means that is not possible to place $n_1$ nodes consistently with the fixed values and constraints imposed by the network topology.

Beside such squares, the phase diagram has some other meaningful areas to discuss. In particular, the high degeneracy squares are considered as the most typical configurations for a random distribution of a property $D = H = 1$; and the phase boundaries squares map atypical configurations. For such phase boundaries different layouts are recognizable. Indeed, in Figure~\ref{fig_config}, four possible configurations are represented (where the point of each configuration is correspondingly numbered in the phase diagram of Figure~\ref{pdexample}): $D \gg 1$ is a dyadic case where black nodes concentrate in a central cluster of the graph which maximizes $m_{11}$; $D \ll1$ is an anti-dyadic configuration where black nodes tend to be farther apart; $H \ll 1$ is an heterophobic configuration where black nodes are located in the peripheral area which minimize $m_{10}$; and, $H \gg 1$ is an heterophilic configuration where black nodes correspond to the most connected nodes so that the edges with white nodes are maximized.

The graphical nature of the phase diagram allows for easier observation of the distribution of the nodes' characteristics, however, as the number of the possible configurations increases exponentially with $N$, the phase diagram is hard to compute for large networks. 
Moreover, since $n_1 \in [0,N]$ it can change while the network structure remains the same. Indeed the number of nodes showing a certain characteristic can change, or different characteristics can be studied by varying $n_1$, such as in~\cite{HAKM}.
Therefore, a complete analysis may require a sequence of $N+1$ phase diagrams computed for each value of $n_1$ as in~\cite{CFI}. An outcome of the latter case for a graph with 25 nodes and 32 edges is shown in Figure~\ref{pdsequence}.

\section{Upper and lower bounds}
\label{ULB}

In this section, we propose an extension of upper bounds (\ref{ub11old}) and (\ref{ub10old}) for a graph $G$ given $N$, $M$ and $n_1 \in [0,N]$. Moreover, we propose  two lower bounds (\textit{LB}) to $m_{11}$ and $m_{10}$. In other words, we want to restrict the feasible region of the dyadic effect for a graph $G$ as much as possible, excluding non-admissible configurations.

We can notice that in~\cite{PB} no lower bounds are provided, thus both bounds on $m_{11}$ and $m_{10}$ are assumed to be zero at minimum. 

\subsection{Upper bound $\UBu$}\label{UBunouno}

Equation (\ref{ub11old}) states that the maximum number of $m_{11}$ within a network is equal to the minimum between two quantities: 
the number of network's edges (meaning that all the edges are $m_{11}$ and $n_1=N$); the number of edges within a clique $K_{n_1}$, i.e. a complete subgraph with $n_1$ nodes within $G$.

The rationale behind the latter statement is that the upper bound is pushed to the maximum value when $G$ is supposed to contain $n_1$ nodes arranged in a clique.

\begin{proposition}\label{propubm11}
Let us consider the degree sequence $D_G$ of the graph $G$ with $n_1$ nodes having the characteristic equal to 1. The upper bound $UBm_{11}$ on the number of edges $m_{11}$ is:
\begin{equation}\label{for_m11}
UBm_{_{11}} = min\Bigg( M, \binom { n_1 } {2}, \bigg\lceil \sum_{i \in \DHn}{\frac {\min(d_{i}, n_1 - 1)}{2}} \bigg\rceil \Bigg)     
\end{equation}
\end{proposition}
\begin{proof}
Considering the degree sequence, we distinguish two different cases based on the fact that $D_G$  may or may not contain $n_1$ nodes of degree at least $n_1 - 1$. In the former case, the necessary condition for the existence of $K_{n_1}$ holds and we can suppose, as a worst case, the realization of a clique $K_{n_1}$ considering $\DHn$. Thus, the bound given in formula (\ref{ub11old}) is the tightest. In the other case, the necessary condition is not upheld, thus we can affirm that $G$ does not contain a clique $K_{n_1}$.

Note that if in $D_G$ it is possible to realize a clique $K_{n_1}$, its nodes can be considered the same as those in $\DHn$. Indeed, $D_G$ is an ordered sequence and if $\DHn$ does not realize $K_{n_1}$, none other subsequence into $D_G$ can realize it. Therefore, we refer our analysis to the subsequence $\DHn$.

As $K_{n_1}$ is the densest possible realization in $\DHn$, we can similarly search for the densest possible realization actually feasible in the subsequence when the clique is not realizable. 

Considering $\DHn$ we distinguish two situations: the sum of its elements being odd or being even. By construction, $D_G$ is graphic but a subsequence can not be graphic albeit the sum of its elements is even. Indeed, $\DHn$ is a graphic sequence only if conditions reported in~\cite{EG, H1962, H1955} hold; otherwise, it is not graphic.

Therefore, let consider the case in which $\DHn$ is graphic. In this situation, the densest hypothetical realization of  $\DHn$ is a graph in which the handshaking lemma holds, thus the number of its edges is $m_{11} = \sum_{i \in \DHn} \min(d_i, n_1 - 1 )/2$ because each one of the $n_1$ nodes has its degree bounded by the value $n_1 - 1$.

Since we are searching for upper bounds, the situation in which $\DHn$ is not graphic can be managed through an overestimation of the handshaking lemma. In this case, we consider the densest hypothetical realization of a simple graph that involves the maximum number of stubs corresponding to $\min(d_i, n_1 -1)$ for any $i \in \DHn$. In order to obtain the $UB$ value, we consider the involvement of all the elements in $\DHn$, ceiling the sum if odd. Through this procedure the handshaking lemma holds and we can compute the number of edges as $m_{11} = \lceil \sum_{i \in \DHn} \min(d_i, n_1 -1) /2\rceil$. 

Summarizing all the considerations thus far, $UBm_{11}$ can be written as in Formula (\ref{for_m11}).
\end{proof}

\subsection{Upper bound $\UBz$}

Equation (\ref{ub10old}) states that the maximum number of $m_{10}$ is equal to the minimum between $M$, meaning that all the edges are $m_{10}$ and thus there are no adjacent $n_1$, and the number of edges within a set of $n_1$ stars of degree $n_0$ (or $n_0$ stars of degree $n_1$).

In more detail, the second element in the upper bound's formula implies that all the $n_1$ nodes are arranged in order to be the central nodes of a set of stars $S_{n_0 + 1}$ with non-adjacent central vertices and of degree $n_0$, or viceversa.

\begin{proposition}
Let us consider the degree sequence $D_G$ of the graph $G$ with $n_1$ nodes having the characteristic equal to 1. The upper bound $UBm_{10}$ on the number of edges $m_{10}$ is:
\begin{equation}\label{for_m10}
UBm_{10} = min\Bigg( M, n_1 n_0, min\bigg( \sum_{i \in D_G^H(n_1)}{\min(d_{i}, n_0)} , \sum_{i \in D_G^H(n_0)}{\min(d_{i}, n_1)}  \bigg)\Bigg)     
\end{equation}
\end{proposition}
\begin{proof}
If the graph $G$ can contain $n_1$ stars $S_{n_0 + 1}$, i.e. in $D_{G}$ are present at least $n_1$ elements with $d_{i} \geq n_0$ then the maximum number of stars is theoretically allowed and the bound given in~\cite{PB} can be considered the tightest. Otherwise if such stars do not exist we can take into account the set of stars that is actually realizable using the degree sequence of $G$. 

Clearly, the same reasoning can be applied when considering $n_0$ instead of $n_1$. For any fixed $n_1$, stars can be realized with all the central nodes having the characteristic $c_i = 1$ and the other elements having $c_i = 0$ and viceversa. 

When maximising $m_{10}$, we ask for the set of stars with non adjacent central nodes that brings $m_{11}$ to be the minimum, i.e. equal to zero. Under these considerations, we are faced with three different situations: $n_1 < n_0$, $n_1 > n_0$ and $n_1 = n_0$.

Suppose that $n_1 < n_0$. In this case, we can partition the degree sequence of $G$ as $D_G = D_G^H(n_1) \cup D_G^T(n_0)$ or as $D_G = D_G^H(n_0) \cup D_G^T(n_1)$ (see Figure~\ref{deg1}).

When $D_G = D_G^H(n_1) \cup D_G^T(n_0)$, the elements in $D_G^H(n_1)$ show a number of stubs equal to the sum of their degree that, in order to realize edges $m_{10}$, have to find their endpoints in $D_G^T(n_0)$. Three cases are admissible:
\begin{itemize}
\item $\sum_{i \in D_G^H(n_1)} d_i =  \sum_{i \in D_G^T(n_0)} d_i $: all stubs in $D_G^H(n_1)$ have an endpoint in $D_G^T(n_0)$. In this case the realization on $D_G$ has $M = m_{10} = \sum_{i \in D_G^H(n_1)} d_i$.
\item $\sum_{i \in D_G^H(n_1)} d_i <  \sum_{i \in D_G^T(n_0)} d_i $: all stubs in $D_G^H(n_1)$ have an endpoint in $D_G^T(n_0)$ but some stubs in $D_G^T(n_0)$ remain free. In this case the realization on $D_G$ has $m_{10} = \sum_{i \in D_G^H(n_1)} d_i$, while $M$ contains some $m_{00}$.
\item $\sum_{i \in D_G^H(n_1)} d_i >  \sum_{i \in D_G^T(n_0)} d_i $: not all stubs in $D_G^H(n_1)$ have an endpoint in $D_G^T(n_0)$. In this case, the realization on $D_G$ has all stubs in $D_G^T(n_0)$ saturated while some residual stubs in $D_G^H(n_1)$ can create edges between those nodes not involved in stars. In this case, $m_{11}$ may be different to zero and an overestimation is $m_{10} = \sum_{i \in D_G^H(n_1)} d_i$.
\end{itemize}

When $D_G = D_G^H(n_0) \cup D_G^T(n_1)$ three cases can be discussed:
\begin{itemize}
\item $\sum_{i \in D_G^H(n_0)} d_i =  \sum_{i \in D_G^T(n_1)} d_i$: since $n_0$ is greater than $n_1$ and $D_G$ is in non-increasing order, this case is not admissible. 
\item $\sum_{i \in D_G^H(n_0)} d_i <  \sum_{i \in D_G^T(n_1)} d_i$: again, this case is not admissible for the same reason as above.
\item $\sum_{i \in D_G^H(n_0)} d_i >  \sum_{i \in D_G^T(n_1)} d_i$: not all stubs in $D_G^H(n_0)$ have an endpoint in $D_G^T(n_1)$. In this case, the realization on $D_G$ has all stubs in $D_G^T(n_1)$ saturated while some residual stubs in $D_G^H(n_0)$ can create edges that increase the number of $m_{00}$. Thus, an overestimation is $m_{10} = \sum_{i \in D_G^H(n_0)} d_i$.
\end{itemize}

Summarizing, we provided certain overestimations on the number of $m_{10}$ and the minimum among them is our upper bound when $n_1 < n_0$. 

When $n_1 > n_0$, we can partition the degree sequence of $G$ as $D_G = D_G^H(n_0) \cup D_G^T(n_1)$ or as $D_G = D_G^H(n_1) \cup D_G^T(n_0)$  (see Figure~\ref{deg2}) and all the same considerations of above can be repeated, using caution to invert $n_0$ and $n_1$. Again, we obtain certain overestimations on the number of $m_{10}$ and the minimum among them is our upper bound.

When $n_1 = n_0$, then  $D_G = D_G^H(n_1) \cup D_G^T(n_0) = D_G^H(n_0) \cup D_G^T(n_1)$. Considering the first partition, three situation can be discussed:

\begin{itemize}
\item $\sum_{i \in D_G^H(n_1)} d_i =  \sum_{i \in D_G^T(n_0)} d_i $: this happens only when $G$ is regular.
\item $\sum_{i \in D_G^H(n_1)} d_i <  \sum_{i \in D_G^T(n_0)} d_i $: this case is not admissible.
\item $\sum_{i \in D_G^H(n_1)} d_i >  \sum_{i \in D_G^T(n_0)} d_i$: not all stubs in $D_G^H(n_1)$ have an endpoint in $D_G^T(n_0)$. In this case, the realization on $D_G$ has all stubs in $D_G^T(n_0)$ saturated while some residual stubs in $D_G^H(n_1)$ can create edges between nodes not involved in stars. In this case $m_{11}$ can be different to zero and an overestimation is $m_{10} = \sum_{i \in D_G^H(n_1)} d_i$.
\end{itemize}

Such considerations can be repeated for the second partition of $D_G$.
Moreover, knowing the size of the two partitions, we can further bound the introduced quantities and, consequently, $m_{10}$. Indeed, in order to realize $m_{10}$, each element in $D_G^H(n_1)$ or $D_G^T(n_1)$ can be connected at most to $n_0$ others while each element in $D_G^H(n_0)$ or $D_G^T(n_0)$ can be connected at most to $n_1$ others. Thus every $d_i \geq n_0$ in $D_G^H(n_1)$ or $D_G^T(n_1)$ is actually bounded by $n_0$ while every $d_i\geq n_1$ in $D_G^H(n_0)$ or $D_G^T(n_0)$ is actually bounded by $n_1$; furthermore the residual degree of each $d_i$ does not contribute to the formation of $m_{10}$. 

Finally, the value for the upper bound on the number of $m_{10}$ can be written as in Formula (\ref{for_m10}).

\end{proof}

\subsection{Lower bound $LBm_{11}$}
\label{LBmu}
We propose a lower bound of $m_{11}$, observing under which conditions a hypothetical graph realization of $D_G$ exists that contains at least some $m_{11}$. Such quantity is considered as an underestimation of $m_{11}$ in the original $G$.

\begin{proposition}
Let us consider the degree sequence $D_G$ of the graph $G$ with $n_1$ nodes having the characteristic equal to 1. The lower bound $LBm_{11}$ on the number of edges $m_{11}$ is:
\begin{equation}\label{lbm11}
 LBm_{11} = \max \Bigg(0, \bigg \lfloor \frac{\sum_{i \in D_G^T(n_1)} d_i -  \sum_{i \in D_G^H(n_0)} d_i }{2} \bigg \rfloor\Bigg)
\end{equation}
\end{proposition}

\begin{proof}
Given $n_1$, let consider again the two possible partitions of the degree sequence $D_G = D_G^H(n_1) \cup D_G^T(n_0) = D_G^H(n_0) \cup D_G^T(n_1)$. The following three cases hold for $n_1 \gtreqless n_0$:

\begin{itemize}
\item $\sum_{i \in D_G^H(n_1)} d_i =  \sum_{i \in D_G^T(n_0)} d_i$ or $\sum_{i \in D_G^H(n_0)} d_i =  \sum_{i \in D_G^T(n_1)} d_i$: when admissible, stubs have endpoints in different partitions and $M = m_{10}$; thus, $m_{11} = 0$. 
\item $\sum_{i \in D_G^H(n_1)} d_i <  \sum_{i \in D_G^T(n_0)} d_i$ or $\sum_{i \in D_G^H(n_0)} d_i > \sum_{i \in D_G^T(n_1)} d_i$: when admissible, there is no room for residual degree in the partition of nodes with the characteristic $c_i = 1$; thus, $m_{11} = 0$.
\item $\sum_{i \in D_G^H(n_1)} d_i >  \sum_{i \in D_G^T(n_0)} d_i$ or $\sum_{i \in D_G^H(n_0)} d_i < \sum_{i \in D_G^T(n_1)} d_i$: when admissible, some stubs in the partition of nodes with the characteristic $c_i = 1$ can link among themselves. 
\end{itemize}

The third case can happen when $\sum_{i \in D_G^H(n_1)} d_i -  \sum_{i \in D_G^T(n_0)} d_i > 0$ or $\sum_{i \in D_G^H(n_0)} d_i - \sum_{i \in D_G^T(n_1)} d_i > 0$.
Note that, since we are searching for the minimum number of stubs able to make $m_{11} \neq 0$ and that for any given $n_1$ $\sum_{i \in D_G^H(n_1)} d_i \geq \sum_{i \in D_G^T(n_1)} d_i$, we can restrict our analysis to the cases when $ D_G = D_G^H(n_0) \cup D_G^T(n_1)$. 
Therefore, a lower bound on the number of $m_{11}$ is given by selecting the maximum value between 0 and an underestimation of the possible edges created in the partition $D_G^T(n_1)$, such as in Formula (\ref{lbm11}). 

\end{proof}
 
\subsection{Lower bound $LBm_{10}$}
\label{LBmz}
When $n_1 = 0$ or $n_1 = N$, the number of $m_{10}$ is trivially 0. In all other cases, any possible connected realization of $D_G$ contains at least an edge with endpoints with different characteristics, i.e. $m_{10} \geq 1$.

In order to discuss a lower bound on $m_{10}$ that, in some cases, overcomes the given inequality from above, we take into account certain arguments based on the realizability of a complete, or at least densest, subgraph from $D_G$, similarly to as in Section~\ref{UBunouno}.

\begin{proposition}
Let us consider the degree sequence $D_G$ of the graph $G$ with $n_1$ nodes having the characteristic equal to 1. The lower bound $LBm_{10}$ on the number of edges $m_{10}$ is:
\begin{equation}
LBm_{10} =
  \begin{cases}
    0                                                                                & \quad \text{if } n_1 = 0, N\\ 
    max \bigg(1; \sum_{i \in D_G^T(n_1)} d_i - n_1(n_1 -1)\bigg)  & \quad \text{if } n_1 \in (0, N)  \\
  \end{cases}
\end{equation}
\end{proposition}
\begin{proof}
When $n_1 \in (0, N)$ we distinguish two cases determined by whether $D_G$ contains $n_1$ nodes of degree at least $n_1 - 1$ that allows for a realization of a clique $K_{n_1}$.

If $D_G$ does not contain at least $n_1$ nodes of degree $n_1 - 1$, we can suppose that the densest hypothetical realization of a simple graph involving the maximum number of stubs is the same as in the proof of Proposition~\ref{propubm11}. Such realization contains the maximum number of $m_{11}$ and at least one residual stub has to find its endpoint in one of the nodes of the remaining part of $D_G$ because any realization should be connected. Therefore, this case leads us to consider $m_{10} = 1$.

If $D_G$ contains at least $n_1$ nodes of degree $n_1 - 1$, then it can hypothetically realize $K_{n_1}$. In order to search for a lower bound on $m_{10}$ we take into account the subsequence which contains the minimum number of stubs, i.e. $D_G^T(n_1)$. 

If a clique can be realized within $D_G^T(n_1)$, then $m_{10} = \sum_{i \in D_G^T(n_1)} d_i - n_1(n_1 - 1)$. Indeed, every stub that constitutes the residual degree will find its endpoint in a node with $c_{i} = 0$. Otherwise, $\sum_{i \in D_G^T(n_1)} d_i < n_1(n_1 -1)$ and $m_{10} = 1$.

\end{proof}

\subsection{Bounds' implications on dyadic effect and its applications}

In Section~\ref{phasdiag} we introduced range values for $D$ and $H$ and their relationships with the corresponding upper bounds. Since we defined two new formulas for $UBm_{11}$ and $UBm_{10}$, $D_{max}$ and $H_{max}$ may assume lower values. Moreover, the introduction of $LBm_{11}$ and $LBm_{10}$ results in the definiton of two new quantities, i.e. $D_{min} = LBm_{11}/ \overline{m}_{11}$ and $H_{min} = LBm_{10}/\overline{m}_{10}$. Thus, we can state that:

\begin{proposition}
Given a simple graph $G$ with $n_1$ nodes having the characteristic equal to 1, the dyadic effect is bounded as follow: $D_{min} \leq D \leq D_{max}$ and $H_{min} \leq H \leq H_{max}$.
\end{proposition}

This proposition has a main implication. Indeed, since the dyadic effect has been used to quantify homophily, the proposition sets the bounds on nodes' tendency to connect with others similar to themselves by using information on the network itself instead of a priori combinatorial arguments. In fact, the bounds presented in~\cite{PB} are valid for every graph of $N$ nodes, $M$ edges and a fixed integer $n_1$, while the bounds presented in Section~\ref{ULB} depend on the graphic sequence of the given graph, thus are valid for the set of all graphs having the same $D_G$. Such set is still wide but has a tighter relationship with the graph under observation.

Regarding the applications, when a network and a set of characteristics are given, it is straightforward to compute for each characteristic the point of maximum degeneration $D = H = 1$, the values $D$, $H$ and through the use of the four bounds, $D_{min}$, $D_{max}$, $H_{min}$ and $H_{max}$. This approach can be useful in many applications, such as in~\cite{HAKM},~\cite{JDC},~\cite{JLJZ} and~\cite{ZZJ} where the phase diagram is hard to compute. In such contexts statistical approaches are used to gather information on the correlation between nodes' characteristics and the network topology by looking for the relative distance of the point $(D, H)$ from the point of maximum degeneration.  

The introduction of the new bounds makes any comparison within the two dimensional space, defined originally in~\cite{PB} and improved in the previous section, more reliable as they are performed on measures that are deeply related to the structure of the analysed graph.

\section{Empirical evidence} 
\label{EE}

Herein, we show empirical evidence computing upper and lower bounds as presented in Section~\ref{ULB} for different networks. In particular, we extensively study the test graph given in Figure~\ref{example} in order to provide an evaluation of the feasible region reductions, then we provide results on different instances in order to observe the behaviour of different bounds.

The data processing, network analysis and all simulations were conducted using the software R~\cite{Rsoft} with the igraph package~\cite{igraph}.

\begin{figure}[htbp]
\begin{center}
\includegraphics[scale=.4]{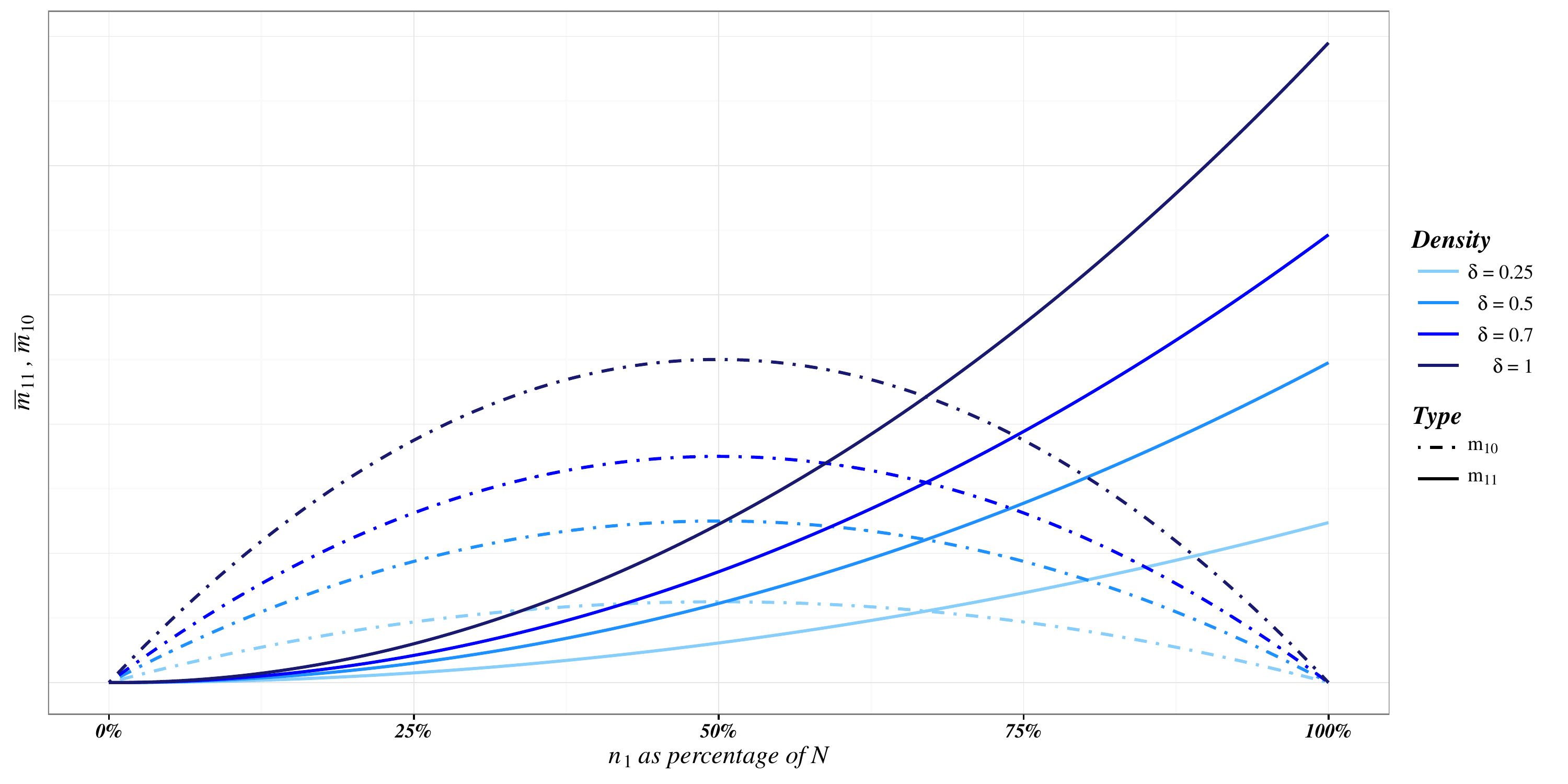}
\caption{$\overline{m}_{11}$ and $\overline{m}_{10}$ values as a function of the fraction of $n_1$ nodes on $N$.}
\label{m10m11}
\end{center}
\end{figure}

All analysis has been conducted considering the given graph $G$ with $N$ nodes and $M$ edges and the value of $n_1$ ranging from $0$ to $N$. Values of $\overline{m}_{11}$ and $\overline{m}_{10}$ were computed straightforwardly and independently of $N$ by means of formulas (\ref{m11}) and (\ref{m10}) and Figure~\ref{m10m11} shows their values as a function of the fraction of $n_1$ nodes on $N$ and for different values of the density $\delta$.

\begin{figure}[htbp]
\begin{center}
 \begin{minipage}[b]{14cm}
   \centering
   \includegraphics[trim=4cm 1.5cm 4cm 1cm, clip=true, totalheight=0.3\textheight]{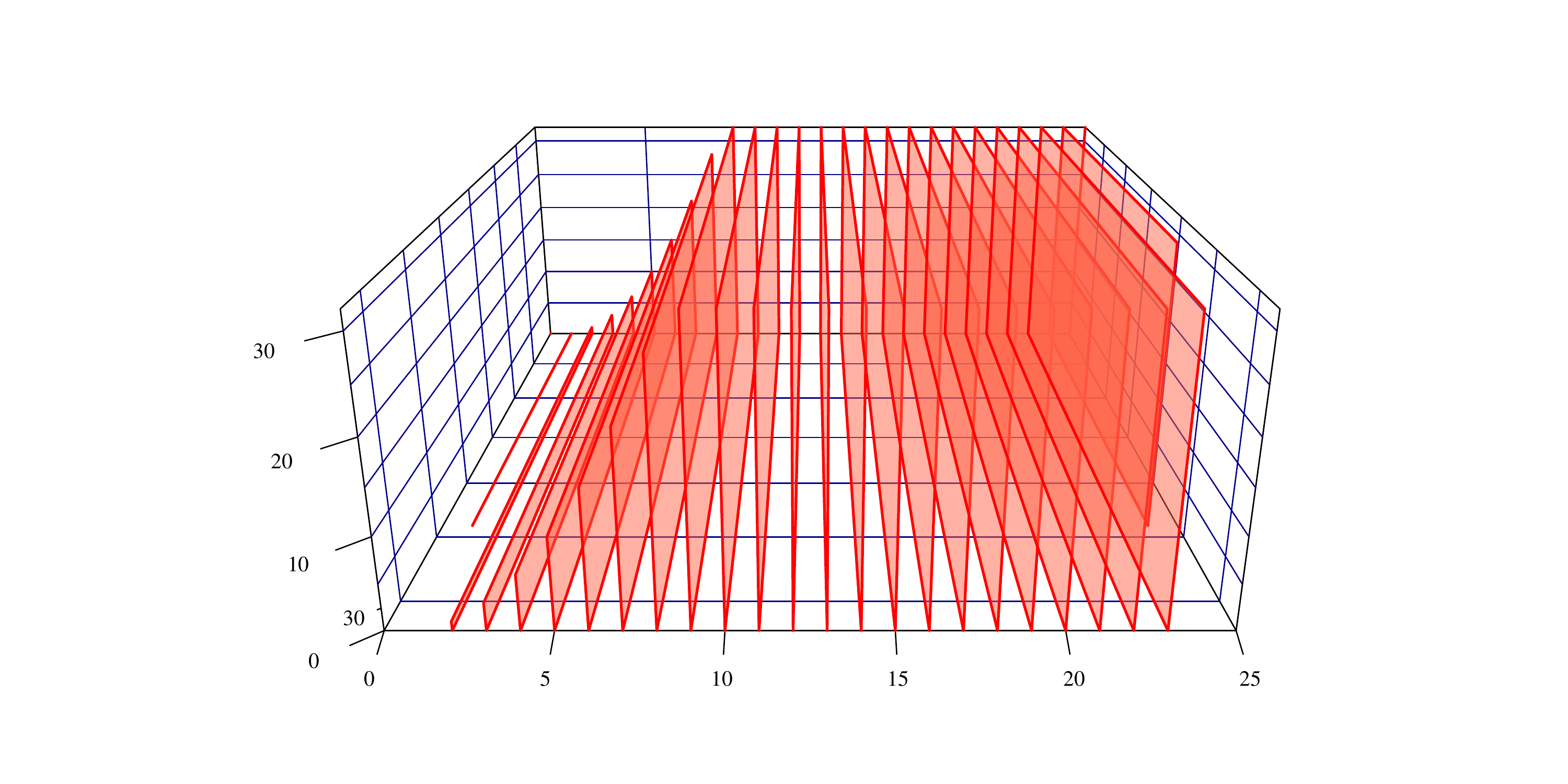}
   \caption{Areas computed with old bounds as a function of $n_1$ nodes on $N$. The $y$ and $z$ axes represent $m_{11}$ and $m_{10}$ respectively.}\label{figgain}
    \label{b_old}
 \end{minipage}
 \\ 
 \begin{minipage}[b]{14cm}
  \centering
   \includegraphics[trim=4cm 1.5cm 4cm 1cm, clip=true, totalheight=0.3\textheight]{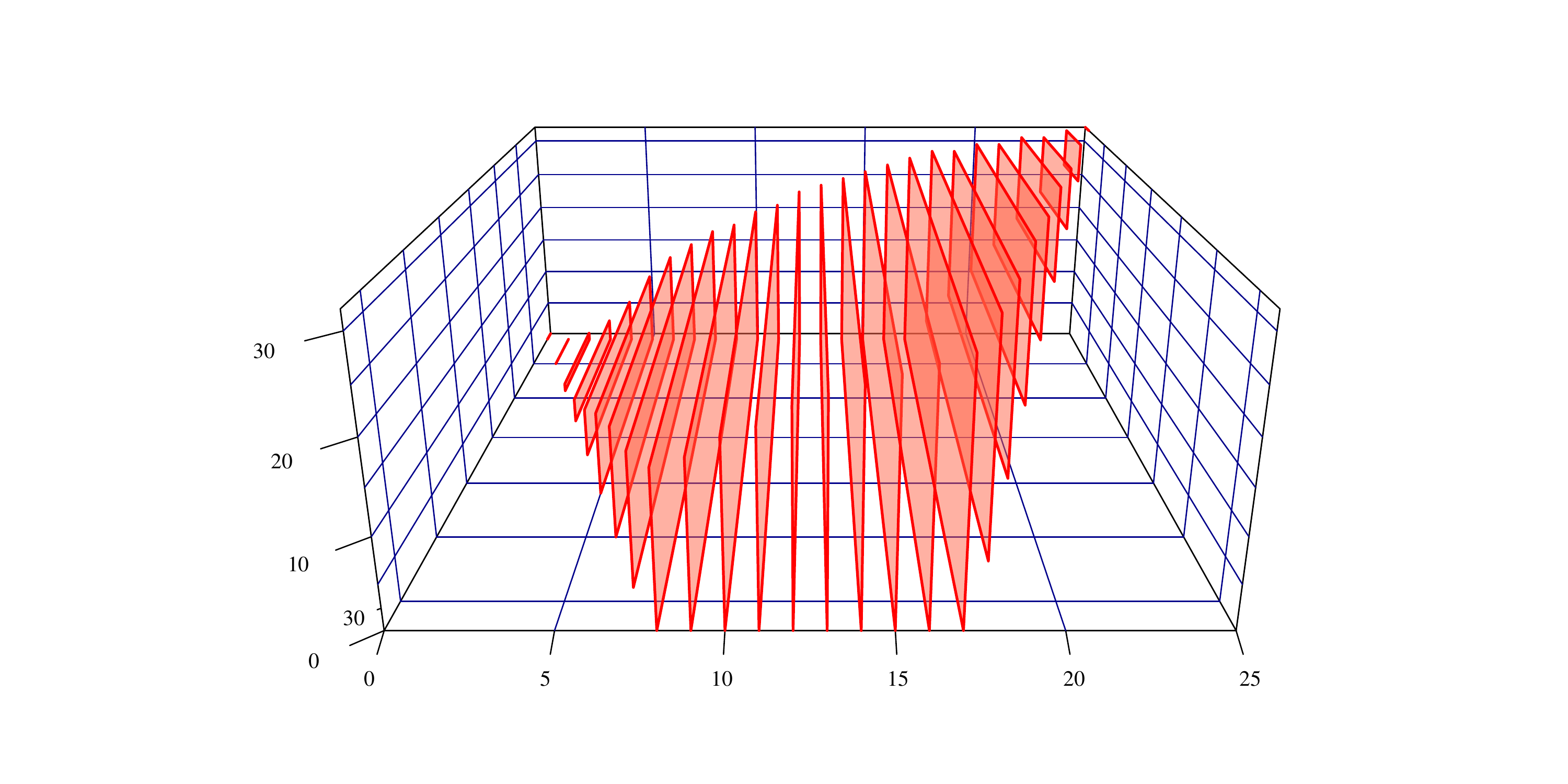}
   \caption{Areas computed with new bounds as a function of $n_1$ nodes on $N$. The $y$ and $z$ axes represent $m_{11}$ and $m_{10}$ respectively.}\label{figold}
   \label{b_new}
 \end{minipage}
\end{center}
\end{figure}

\subsection{Analysis of the test graph}
We begin by illustrating the shrinkage of the area in which the phase diagram lies by using the test graph with $N = 25$ and $M = 32$. The areas in Figure~\ref{b_old} and~\ref{b_new} represent the sequence of feasible regions for $n_1$ varying from $0$ to $N$, each bounded by the correspondent values of $UBm_{10}$, $LBm_{10}$, $UBm_{11}$ and $LBm_{11}$. Comparing Figure~\ref{b_old} and Figure~\ref{b_new} we can immediately notice the difference between the feasible regions provided by old and new bounds as well as the consequent improvements mostly appreciable for high and low values of $n_1$. In Figure~\ref{b_new} we observe how, by applying the bounds, the feasible region changes together with the various shapes of the $N+1$ phase diagrams as shown in Figure~\ref{pdsequence}. Moreover, the areas evolve following a trajectory which reflects the trend observed in the curves of the expected values of $m_{11}$ and $m_{10}$ in  Figure~\ref{m10m11}.       

\subsection{Analysis of benchmark instances}

We applied the proposed bounds to three classes of networks of various size according to their degree distribution. We chose Erd\H{o}s and R\'enyi random graphs~\cite{ER} and scale-free networks~\cite{BA} and, for each class, we generated ten instances and the results herein are presented as averaged. We considered networks with $N = 1000, 3000$ with two different settings: the first with an average degree $\langle d \rangle = 6$; the second having a density $\delta = 0.9$.  In addition, for each setting we included regular graphs. The first setting was considered in order to perform an analysis similar to~\cite{PB} while the second setting was chosen in order to test bounds that require $D_G$ from dense graphs.

\begin{figure}[htbp]
\begin{center}
 \begin{minipage}[b]{14cm}
   \centering
   \includegraphics[scale=.45]{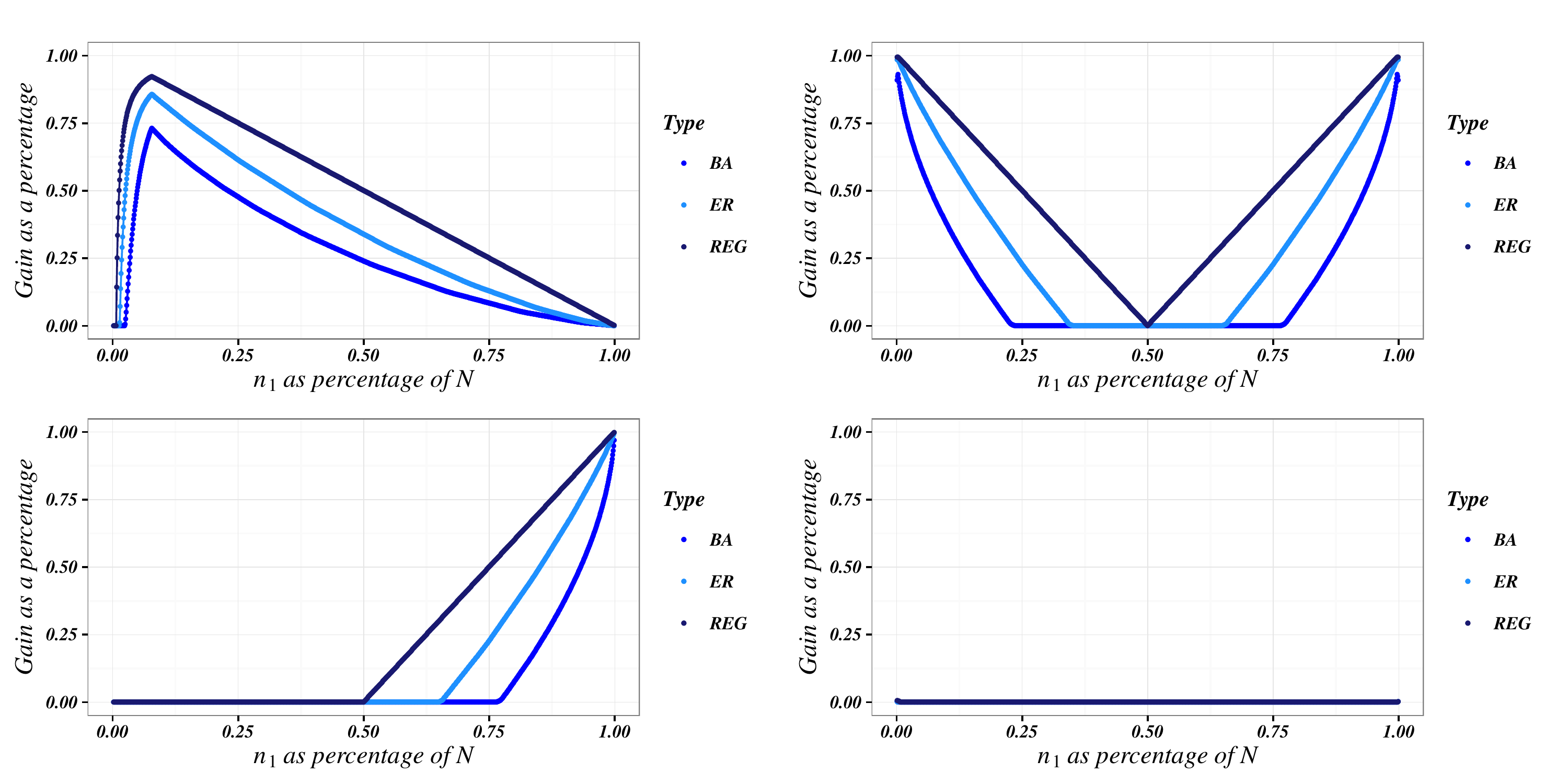}
   \caption{Clockwise representation of the gain deriving from the application of $UBm_{11}$, $UBm_{10}$, $LBm_{11}$, $LBm_{10}$ as a function of $n_1$ nodes on $N = 1000$. Average degree was set to $\langle d \rangle = 6$.}
    \label{grid_1000}
 \end{minipage}
 \begin{minipage}[b]{14cm}
  \centering
   \includegraphics[scale=.45]{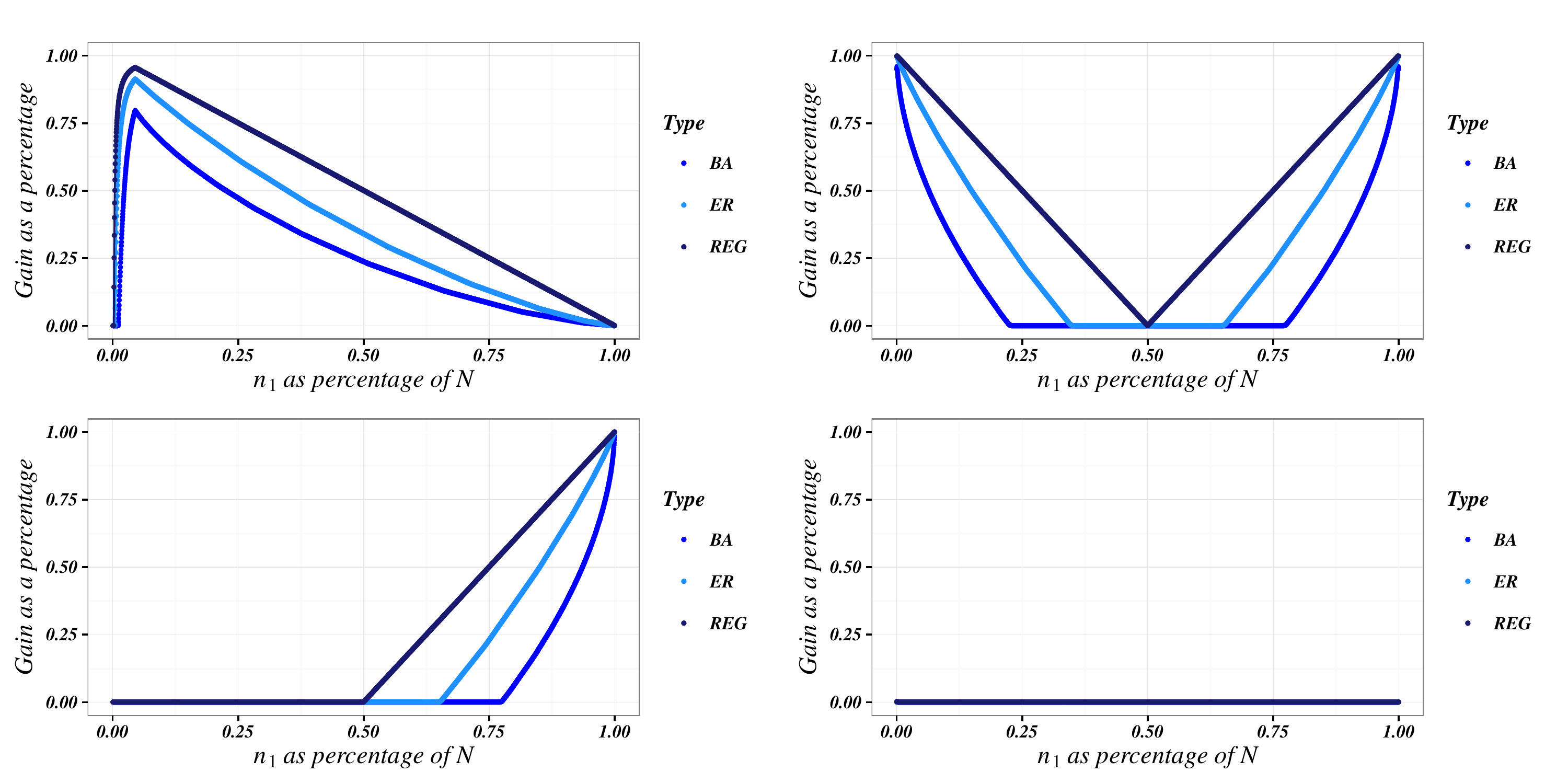}
   \caption{Clockwise representation of the gain deriving from the application of $UBm_{11}$, $UBm_{10}$, $LBm_{11}$, $LBm_{10}$ as a function of $n_1$ nodes on $N = 3000$. Average degree was set to $\langle d \rangle = 6$.}
\label{grid_3000}
 \end{minipage}
 \end{center}
\end{figure}

\begin{figure}
\begin{center}
 \begin{minipage}[b]{14cm}
  \centering
   \includegraphics[scale=.45]{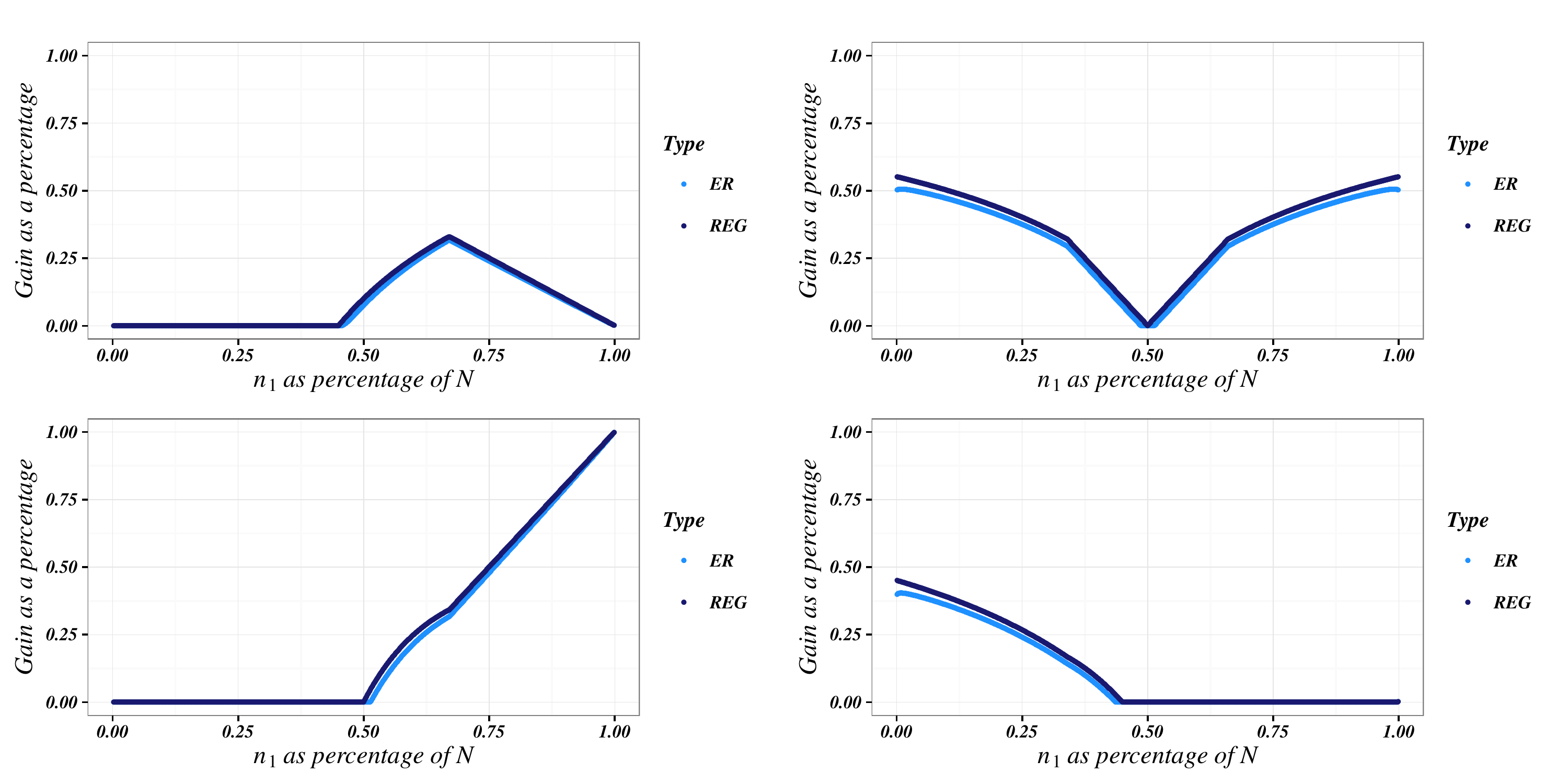}
   \caption{Clockwise representation of the gain deriving from the application of $UBm_{11}$, $UBm_{10}$, $LBm_{11}$, $LBm_{10}$ as a function of $n_1$ nodes on $N = 1000$. Density was set to $\delta = 0.9$.}
\label{grid_1000_d09}
 \end{minipage}
 \begin{minipage}[b]{14cm}
  \centering
   \includegraphics[scale=.45]{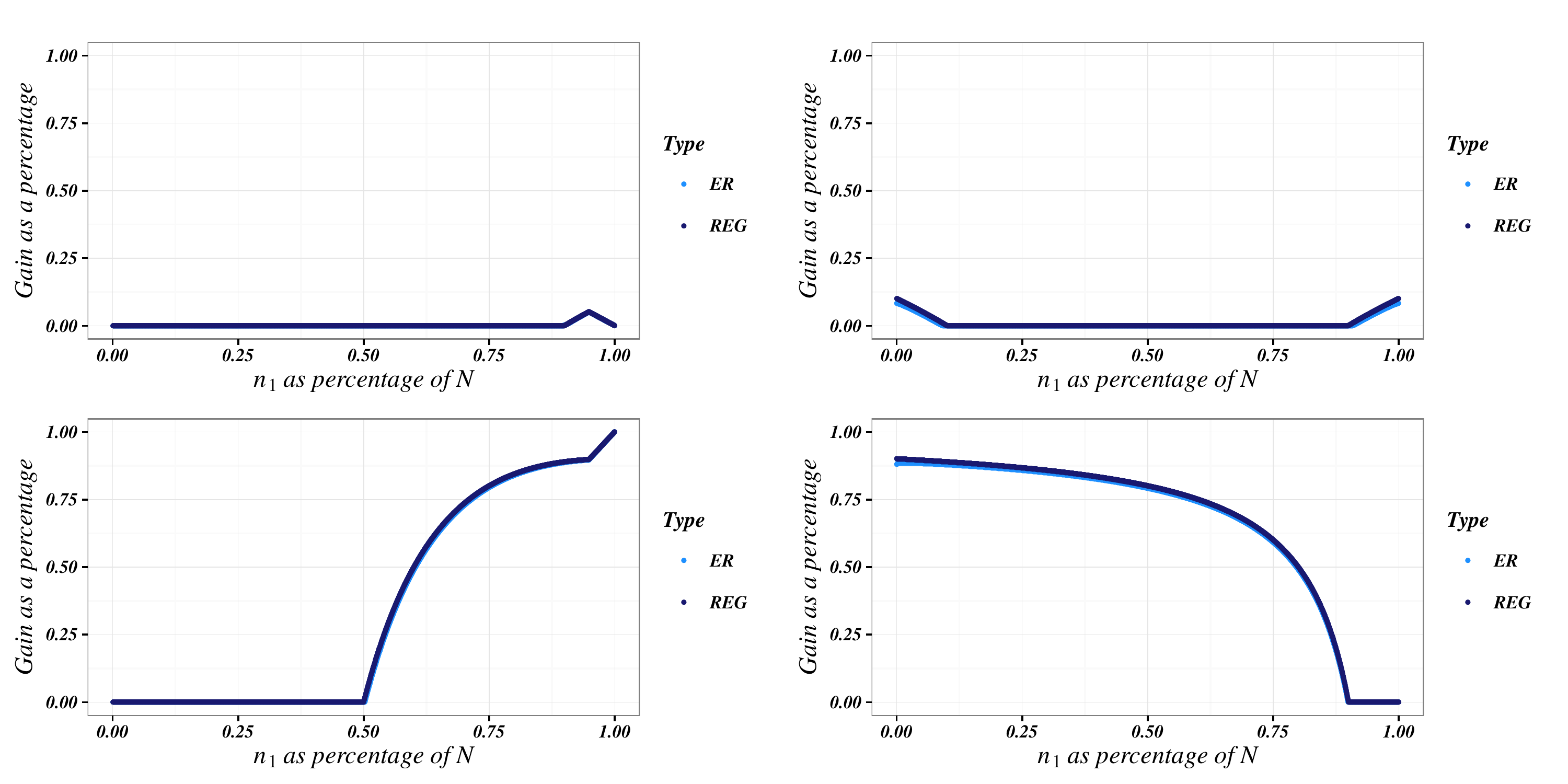}
   \caption{Clockwise representation of the gain deriving from the application of $UBm_{11}$, $UBm_{10}$, $LBm_{11}$, $LBm_{10}$ as a function of $n_1$ nodes on $N = 3000$. Density was set to $\delta = 0.9$}
\label{grid_3000_d09}
 \end{minipage}
\end{center}
\end{figure}

The curves in Figure~\ref{grid_1000} and~\ref{grid_3000} show the percentage gain obtained from applying each new bound from the perspective of the area covered by the feasible region, with respect to the bounds in Formula (\ref{ub11old}) and (\ref{ub10old}) when $\langle d \rangle = 6$. This kind of analysis allows us to estimate the bound's behaviour and the dependency of the gain for different networks types. Indeed, we can firstly observe how, fixed the mean degree, the networks' size do not affect the bounds as behave exactly the same, while the value of $n_1$ acts as a threshold to trigger the bounds. Indeed, upper bounds of Formula (\ref{for_m11}) and (\ref{for_m10}) have a non homogeneous behaviour and for some values of $n_1$ tend to get closer to the values achieved by the bound in Formula (\ref{ub11old}) and (\ref{ub10old}). 

Through observing the analytical relationships referred the upper bounds, it is evident that homogeneous $D_G$ tends to perform better, while those of a more heterogeneous nature such as, for instance, scale-free networks, perform slightly worse due to the deviation of some nodes from the mean degree. 

Lower bound $LBm_{11}$ evaluation can be performed using a similar approach as for the upper bounds. Note that we can not make any comparisons with previous results being implicitly set to $0$ in~\cite{PB}. 

$LBm_{10}$ gives a contribution that can be considered close to zero due to the fact that it is computed as the difference between the tail of $D_{G}$ and the edges within $K_{n_1}$. Indeed, when $\langle d \rangle = 6$ we deal with very sparse graph, thus the latter difference is always negligible.

We tuned the density to high values ($\delta = 0.9$) since the terms in $LBm_{10}$ formula depend on high degree elements within $D_G$. Figures~\ref{grid_1000_d09} and~\ref{grid_3000_d09} show the gain obtained through the introduction of new bounds for both Erd\H{o}s-R\'enyi and regular graphs. Scale-free networks were omitted in the test as they are characterized by low density values~\cite{DGB}.

In this case, the impact of both lower bounds is more significant, particularly for $LBm_{10}$ which differs to zero most substantially when $N$ grows. Indeed, as we noticed before, all the bounds depend on the mean degree  $\langle d \rangle$. Therefore, to keep the value of the density $\delta$ while increasing the network's size, implies the growth of the mean degree as well, since $\langle d \rangle = \delta (N-1)$. Comparing Figures~\ref{grid_1000_d09} and~\ref{grid_3000_d09} it is evident that when the value of $N$ is higher, the impact of the upper bounds tends to be more negligible since their gain is substituted by that of the lower bounds. Finally, the performed analysis reinforces the importance of introducing lower bounds, especially when dealing with dense graphs.

\section{Conclusions}
\label{conclusion} 
In this paper, we developed two upper bounds and two lower bounds in order to reduce the area of the two dimensional space used to represent the feasible region of the dyadic effect in a network with certain nodes' characteristics. Using commonly accepted structural principles, we improved the upper bounds and provided two new lower bounds. The four bounds can be computed using straightforward analytical relationships with no restrictions on either the network size or classes with the only limitation for the graph to be simple. These aspects are particularly relevant in the applications where large networks are investigated taking into account numerous characteristics and  the dyadic effect is studied avoiding the computation of the phase diagram.

The computational analysis of various classes of networks resulted in behavioural differences depending on the inner structure and on the shape of the degree sequence.

Further research should be devoted to study additional aspects. In particular, the asymptotic behaviour of each bound when the network size grows should be considered. Another point of interest would be to improve the bounds using further arguments, especially in the case of scale-free networks where the proposed approach have seemed to perform less effectively.

Finally, considering the strict relationship between assortative mixing and dyadic effect, the proposed bounds and their implications could lay the bases for potential studies in such direction.


\end{document}